\documentclass[11pt]{article}

\usepackage{amsmath,amsfonts,amssymb,amsthm,mathtools,bbm, geometry}

\usepackage[all]{xy}

\def\E{{\cal E}} 
\def\L{{\cal L}}

\usepackage{amssymb}

\newtheorem{thm}{Theorem}[section]

\newtheorem{dfn}[thm]{Definition}

\newtheorem{lem}[thm]{Lemma}

\newtheorem{prop}[thm]{Proposition}

\newtheorem{remark}[thm]{Remark}

\newtheorem{cor}[thm]{Corollary}

\newtheorem{ex}[thm]{Example}

\newtheorem{question}[thm]{Question}

\def\sq{{\scriptscriptstyle \square}}
\def\di{\diamond}

\def\bq{\begin{question}}
\def\bt{\begin{thm}}
\def\bp{\begin{prop}}
\def\blem{\begin{lem}}
\def\bd{\begin{dfn}}
\def\br{\begin{remark}}
\def\bc{\begin{cor}}
\def\bex{\begin{ex}}
\def\beqs{\begin{eqnarray*}}
\def\beq{\begin{eqnarray}}
\def\bi{\begin{itemize}}

\def\eq{\end{question}}
\def\et{\end{thm}}
\def\ep{\end{prop}}
\def\elem{\end{lem}}
\def\ed{\end{dfn}}
\def\er{\end{remark}}
\def\ec{\end{cor}}
\def\eex{\end{ex}}
\def\eeqs{\end{eqnarray*}}
\def\eeq{\end{eqnarray}}
\def\ei{\end{itemize}}

\def\ds{\displaystyle}

\def\c{\cdot}

\def\r{\rangle}
\def\l{\langle}

\def\fB{{\mathfrak B}}

\def\w*{w^*-w^*}

\def\ra{\rightarrow}

\def\vp{\varphi}
 
\def\bs{\backslash}

\def\om{\omega} 
\def\M{{\cal M}}
\def\fC{\mathfrak{C}} 
\def\fS{\mathfrak{S}} 

\def\E{{\cal E}} 
\def\L{{\cal L}}

\makeatletter
\def\moverlay{\mathpalette\mov@rlay}
\def\mov@rlay#1#2{\leavevmode\vtop{%
   \baselineskip\z@skip \lineskiplimit-\maxdimen
   \ialign{\hfil$\m@th#1##$\hfil\cr#2\crcr}}}
\newcommand{\charfusion}[3][\mathord]{
    #1{\ifx#1\mathop\vphantom{#2}\fi
        \mathpalette\mov@rlay{#2\cr#3}
      }
    \ifx#1\mathop\expandafter\displaylimits\fi}
\makeatother

\begin{document}

\title{ On Beurling Measure Algebras }
\author{Ross Stokke\footnote{
This research was partially supported   by an NSERC grant. }    }
\date{}
\maketitle

\begin{abstract}{\small  We show how the measure theory of regular compacted-Borel measures defined on the $\delta$-ring of compacted-Borel subsets of a weighted locally compact group $(G,\omega)$ provides a compatible framework for defining the corresponding Beurling measure algebra $\M(G,\omega)$, thus  filling a gap in the literature.   

\smallskip

\noindent{\em Primary MSC code:} 43A10 \ \ \ {\em Secondary MSC codes:} 22D15, 43A05,  43A20, 43A60, 28C10  \\
{\em Key words and phrases:} Weighted locally compact group, group algebra, measure algebra, Beurling algebra 

 }
\end{abstract}

Throughout this article, $G$ denotes a locally compact group and $\omega: G \ra (0,\infty)$ is a continuous \it weight function \rm  satisfying
$$\om(st) \leq \om(s)\om(t) \ \ (s, t \in G) \qquad {\rm and} \qquad \omega(e_G) = 1;$$
the pair $(G, \om)$ is called a \it weighted locally compact group. \rm   Let $\lambda$ denote a fixed Haar measure on $G$, with respect to which the group algebra $L^1(G)$ and $L^\infty(G)= L^1(G)^*$ are defined in the usual way.  The Beurling group algebra,   $L^1(G, \om)$,  is composed of all functions $f$ such that   $ \omega f$ belongs to $L^1(G)$, with $\|f\|_{1,\om} := \|\om f\|_1$ and convolution product. If   ${\cal S}(G)$   is a closed subspace of $L^\infty(G)$,  $\psi \in {\cal S}(G, \om^{-1})$ exactly when ${\psi \over \om} \in {\cal S}(G)$;  putting $\|\psi \|_{\infty, \om^{-1}} = \left\| {\psi\over \om}\right\|_\infty$,  ${\cal S}(G, \om^{-1}) $ is a Banach space and 
$ S : {\cal S}(G, \om^{-1})  \ra {\cal S}(G): \psi \mapsto {\psi\over \om}$ is an isometric linear isomorphism. 
The Beurling group algebra $L^1(G, \om)$ 
   has become a classical object of study that has received significant research attention over the years:  see the monographs  \cite{Dal-Lau, Kan, Rei-Ste} and the references therein; a sample of relevant articles include  \cite{Gha1,Gha2, Gha-Zad, Gro,  Sam, She-Zha, Zad}.    
  When $\omega$ is the trivial weight $\omega \equiv 1$ ---  the ``non-weighted case" --- $L^1(G, \om) = L^1(G)$, the study of which is intimately linked with the measure algebra $M(G)$ of complex, regular, Borel measures on $G$, which contains $L^1(G)$ as a closed ideal.

     The above definition of $L^1(G, \om)$  is valid for any weight $\omega$.  As in the non-weighted case, it is desirable to have a Beurling measure algebra $M(G,\om)$  that shares the same relationship with  $L^1(G, \om)$ that $M(G)$ shares with  $L^1(G)$.   In the literature, $M(G, \om)$ is usually defined as the collection of all complex regular measures  $\nu$ defined on $\fB(G)$, the $\sigma$-algebra of Borel subsets of $G$,  such that $ \int \om(t) \, d |\nu|(t) < \infty, $ and the identification $M(G, \om) = C_0(G, \om^{-1})^*$  through $\l \nu, \psi\r_\om = \int \psi \, d \nu$ is required.  This implies that the dual map, $S^*$, of the isometric isomorphism $S: C_0(G, \om^{-1}) \ra C_0(G)$ is itself a linear isometric isomorphism of $M(G)$ onto $M(G,\om)$.  Validity of this definition of $M(G, \om)$ thus requires that for each $\mu \in M(G)$, $\nu = S^*\mu \in M(G, \om)$ is a complex  Borel measure defined on all of $\fB(G)$ --- the near-universal requirement of “Borel measures” in abstract harmonic analysis --- satisfying  
  \beq   \label{M(G,w) original def Eqn}  \int \psi \, d \nu = \l \nu, \psi\r_\om = \left\l \mu , {\psi \over \om} \right\r = \int {\psi \over \om} \, d \mu \qquad (\psi \in C_0(G, \om^{-1})).  \eeq 
  However, when $\om$ is not bounded away from zero, it can happen that no such complex measure on $\fB(G)$ exists.   
  
  To see this, consider $(G, \om)$ where $G = (\mathbb{Z}, +)$ and $\om(n) = 2^{-n}$ $(n \in \mathbb{Z})$, and assume the above definition of $M(G, \om)$ is sound.  Since $\mu_1, \mu_2 \in \ell^1(\mathbb{Z})^+ = M(\mathbb{Z})^+$  and $\mu = \mu_1 - \mu_2 \in M(\mathbb{Z})$, where 
  $$\displaystyle{\mu_1(n) = \left \{ \begin{array}{ll}
                          2^{-n}     & \mbox{ $ n \in  2 \mathbb{N}$}\\
                           0 & \mbox{  otherwise}
                            \end{array}
                 \right. }     \qquad {\rm and}  \qquad    \displaystyle{\mu_2(n) = \left \{ \begin{array}{ll}
                          2^{-n}     & \mbox{ $ n \in  \mathbb{N} \bs 2 \mathbb{N}$}\\
                           0 & \mbox{  otherwise}
                            \end{array}
                 \right. },$$
 $\nu_1 = S^*(\mu_1)$,   $\nu_2 = S^*(\mu_2)$, and $\nu = S^*(\mu) = \nu_1 - \nu_2$ are then required to be complex measures on $\fB(G) = \wp(\mathbb{Z})$ satisfying (\ref{M(G,w) original def Eqn}). Hence, for each $n \in \mathbb{Z}$, 
         $$\displaystyle{\nu_1(\{n\}) = \int \chi_{\{n\}} \, d \nu_1 = \left\l \mu_1, {\chi_{\{n\}} \over \om} \right\r =  \left \{ \begin{array}{ll}
                          1    & \mbox{ $ n \in  2 \mathbb{N}$}\\
                           0 & \mbox{  otherwise}
                            \end{array}
                 \right. }  \quad {\rm and}  \quad  \displaystyle{\nu_2(\{n\}) = \left \{ \begin{array}{ll}
                          1     & \mbox{ $ n \in  \mathbb{N} \bs 2 \mathbb{N}$}\\
                           0 & \mbox{  otherwise}
                            \end{array}
                 \right. }; $$
 therefore, $\ds \nu_1(2\mathbb{N}) = \sum_{k \in \mathbb{N}} \nu_1(\{2k\}) = +\infty$ and $\ds \nu_2( \mathbb{N} \bs 2\mathbb{N}) = \sum_{k \in \mathbb{N}} \nu_2(\{2k-1\}) = +\infty$. Thus, $\nu_1$, $\nu_2$ do not map into $\mathbb{C}$. Moreover, (although $\nu_1$, $\nu_2$ can be viewed as positive measures),  if $\nu = \nu_1 - \nu_2$ were a measure, additivity would give $$\nu(\mathbb{N}) = \nu(2\mathbb{N}) + \nu (\mathbb{N} \bs 2\mathbb{N}) = \nu_1(2\mathbb{N}) - \nu_2(\mathbb{N} \bs 2\mathbb{N}) = \infty - \infty.$$  
 
 We conclude that functionals in $C_0(G, \om^{-1})^*$ cannot necessarily be identified with complex Borel measures in the standard sense.   It is perhaps for this reason that many authors assume the additional condition $\om \geq 1$, since this guarantees containment of  $M(G,\om)$ in $M(G)$ and, thus, the essential properties of $M(G)$ also hold for $M(G,\om)$, e.g., see \cite{Dal-Lau}. 
 Letting  $\fS(G)$ denote the $\delta$-ring of ``compacted-Borel sets" --- i.e., the $\delta$-ring of all Borel subsets of $G$ with compact closure --- a \it compacted-Borel measure \rm on $G$ is a countably additive complex-valued function on $\fS(G)$ in the sense of \cite[Definitions II.1.2 and II.8.2]{Fel-Dor}\footnote{In \cite{Fel-Dor}, for the sake of brevity, the authors refer to compacted-Borel measures simply as Borel measures. To our knowledge, with the exception of \cite{Fel-Dor}, Borel measures in abstract harmonic analysis  are always defined on $\mathfrak{B}(G)$.}.  For non-compact $G$, there are positive regular measures $\mu, \nu$ on $\fB(G)$ such that $\mu(G) = \nu(G) = \infty$ (e.g., Haar measures), and therefore $\mu - \nu$ is not defined on $\fB(G)$;  however,  these same measures are real-valued on $\fS(G)$, so  $\mu-\nu$ is well-defined on $\fS(G)$. This is one benefit to studying measure theory over $\fS(G)$, rather than on all of $\fB(G)$.

    The purpose of this article is to show that the theory of complex regular compacted-Borel measures, as developed in \cite{Fel-Dor} (also see paragraph two of the ``Notes and Remarks" section of Chapter II of \cite{Fel-Dor} for additional references), can be used to provide a rigorous definition of $M(G, \omega)$,  thus providing a solid foundation for all the papers in which $M(G, \om)$ is employed without the requirement that $\om \geq 1$; moreover, we hope this reduces the number of instances in which the $\om \geq 1$ assumption is required going forward. To stress that we are using the theory of complex regular compacted-Borel measures,  we will use the notation $\M(G,\om)$ --- inspired by \cite{Fel-Dor} ---   rather than $M(G,\om)$.  Beyond identifying the correct collection of measures to employ, work is required to establish the needed theory. As measure theory can be quite  finicky in general;  because the study of compacted-Borel measures introduces different technicalities than those encountered in the Borel measure situation; and  because a lot of research already depends on  the results found herein, we have included a careful  treatment of our development of $\M(G, \om)$. 
        There are numerous detailed classical expositions of the  basic theory $M(G)$, and we believe the same is required for $\M(G, \om)$. 
  
  We restrict ourselves to developing only the most standard properties of $\M(G, \om)$: we provide a careful definition of its elements and show that with  convolution product it is a dual Banach algebra  containing a copy of the Beurling group algebra $L^1(G, \om)$ as a closed ideal. Beyond this, we only show that $\M(G,\om)$ embeds via a strict-to-weak$^*$ continuous isometric isomorphism as a subalgebra of the universal enveloping dual Banach algebra of $\L^1(G, \om)$, $WAP(L^\infty(G, \om))^*$, a result needed in \cite{Kro-Ste-Sto-Yee}. The inspiration for this paper was our need to work with $\M(G, \om)$ in \cite{Kro-Ste-Sto-Yee}. 
  
  
  \section{$\M(G, \om)$: definition and basic properties }

     Unless explicitly indicated otherwise, \it  all references are to statements  in  \S s 1,2,5,7-10 of Chapter II and \S 10 of Chapter III  of \rm \cite{Fel-Dor}.  We will mostly   adhere to the notation found therein.  In particular, $\M(G)$ is the linear space composed of all regular complex compacted-Borel measures on $G$ (\S s II.8 and III.10) and $\M_r(G)$ is the Banach space of \it bounded \rm measures in $\M(G)$ (\S s II.1 and II.8).  Let $\fC(G)$ denote the directed set of compact subsets of $G$, and denote  the space of continuous functions on $G$ with compact support by $C_{00}(G)$,   the space of continuous functions on $G$ vanishing at infinity by $C_0(G)$, and the space of continuous functions  on $G$ supported on $K \in \fC(G)$ by $C_K(G)$; unless the context requires otherwise,  these spaces are  taken with the uniform norm $\|\c \|_\infty$.

 \br   \label{Max regular extension remark}   \rm (a) Let $\mu \in \M(G)$.   A Borel subset $A$ of $G$ belongs to ${\cal E}_\mu$ if $A$ is contained in some open set $U$ such that 
 $$\sup \{ |\mu|(A'): A' \in \fS(G) \ {\rm and } \ A' \subseteq U \} < \infty;$$
  ${\cal E}_\mu$ is a $\delta$-ring containing $\fS(G)$ and,   for $A \in {\cal E}_\mu$, putting 
 \beq \label{Max regular extension eqn} \mu_e(A):= \lim_C \mu(C) \quad {\rm where} \ C \in \fC(G), \ C \subseteq A,\eeq
 we obtain a complex measure on ${\cal E}_\mu$ extending $\mu$, called the maximal regular extension of $\mu$ (II.8.15).  Observe that any Borel subset of a set in ${\cal E}_\mu$ is also in ${\cal E}_\mu$, from which it readily follows that $h \chi_E$ is locally $\mu_e$-measurable whenever  $E \in {\cal E}_\mu$ and $h$ is a Borel-measurable function on $G$. 
  
 \smallskip
 
 \noindent (b)   When $\mu \in \M_r(G)$, ${\cal E}_\mu = \fB(G)$ and $\mu_e \in M(G)$, where $M(G)$ denotes the usual measure algebra of regular complex Borel measures $\mu: \fB(G) \ra \mathbb{C}$, e.g., see \cite{Dal, HR, Pal}. Thus, the measures  in  $\M_r(G)$ are in one-to-one correspondence with measures  in $M(G)$ via $\mu \mapsto \mu_e$; moreover, it is clear from the results in \S III.10 (or Theorem \ref{M(G,w) Banach Algebra Thm}, below,  in the non-weighted case) that $\mu \mapsto \mu_e$ is a weak$^*$-continuous isometric algebra isomorphism of $\M_r(G)$ onto $M(G)$.  Thus, for the purposes of abstract harmonic analysis on (non-weighted) $G$, $\M_r(G)$ can be used in place of the usual $M(G)$, and, as shown in \cite{Fel-Dor},  provides some advantages. 
 
 \er

For $\mu \in \M(G)$, let $I_\mu$ denote the linear functional $ I_\mu(f) = \int f \, d\mu$ defined on $\L^1(\mu)$, or any subspace of $\L^1(\mu)$.  Then \beq \label{Riesz Rep Eqn} \mu \mapsto I_\mu : \M(G) \ra {\mathfrak I} \eeq
is a linear bijection where $\mathfrak I$ is the set of all linear functionals $I$ on $C_{00}(G)$ such that $I \in C_K(G)^*$ for each $K \in \fC(G)$;  (\ref{Riesz Rep Eqn}) maps $\M(G)^+$ onto ${\mathfrak I}^+$ and $\M_r(G)$ onto $C_{00}(G)^* = C_0(G)^*$  (II.8.12).  

\br \rm    \label{M_r=C0* Remark}   It should be noted that when $\mu$ is a complex measure on a $\delta$-ring $\fS$, $f \in \L^1(\mu)$ requires that $f$ vanish off a countable union of sets in $\fS$ (II.2.5, paragraph 2). Thus, when $f \in \L^1(\mu)$ for $\mu \in \M(G)$, $f$ must vanish off a $\sigma$-compact set, a technical issue  requiring  careful attention throughout this note.   Consider the case when  $\mu \in \M_r(G)$. Then any $\phi \in C_0(G)$ vanishes off a $\sigma$-compact set and since $\phi$ is continuous and  bounded, it is easy to see that $\phi \in \L^1(\mu)$. Assuming further that $\mu \geq 0$ and $\phi \geq 0$ and taking an increasing sequence $(\phi_n)$ in $C_{00}(G)^+$ such that $\| \phi_n - \phi \|_\infty \ra 0$, $\lim I_\mu(\phi_n) = \lim \int \phi_n \, d \mu = \int \phi \, d \mu = I_\mu (\phi)$ (e.g., by MCT II.7), so $I_\mu$ is the unique continuous extension of $I_\mu$ on $C_{00}(G)$ to $C_0(G)$. Thus, $C_0(G)^* = \{I_\mu: \mu \in \M_r(G) \}$, so --- in this theory and as usual --- we can identify $\M_r(G)$ and $C_0(G)^*$ through the pairing $\l \mu, \phi \r = \int \phi \, d \mu$.  
 \er 

Let $\nu \in \M(G)$, $h$ a continuous function on $G$.   Then $h$ is locally $\nu$-measurable (II.8.2) and for each $A \in \fS(G)$, $h \chi_A \in \L^1(\nu)$ since $|h|$ is bounded on $A$; i.e., $h$ is locally $\nu$-summable. Therefore, $$h \nu(A): = \int h \chi_A \, d \nu \qquad (A \in \fS(G))$$ defines a complex measure on $\fS(G)$ (see II.7.2, where the notation $h \, d\nu$ rather than $h \nu$ is used); as $h \nu \ll \nu$ (II.7.8), $h \nu \in \M(G)$ (II.8.3). If $h >0$, then $ {1 \over h}(h \nu) \in \M(G)$ and a simple application of II.7.5 gives  ${1 \over h}(h \nu) = \nu$.  

Hence, $\om \nu \sim \nu$
for each $\nu \in \M(G)$, and 
$$\M(G) \ra \M(G): \nu \mapsto \om \nu$$ 
defines a linear isomorphism with inverse $\nu \mapsto {1 \over \om} \nu$.  We can thus define
$$\M(G,\om):= \{ \nu \in \M(G): \om \nu \in \M_r(G)\}; \quad   {\rm  letting}  \ \ \  \| \nu\|_\om = \|\om \nu \| \ \ \  (\nu \in \M(G,\om)),$$ it follows that $\M(G,\om)$ is a Banach space and $\nu \mapsto \om \nu$ is an isometric linear isomorphism of $\M(G, \om)$ onto $\M_r(G)$ with inverse map $\mu \mapsto {1 \over \om} \mu$.  
(As shown in the introduction, this definition cannot, in general, be made  with $M(G)$ replacing $\M_r(G)$.) Observe that by II.7.3,  $\nu \in \M(G,\om)$ exactly when $|\nu | \in \M(G,\om)$, and $\| \nu \|_\om = \| | \nu | \|_\om$.

\bp \rm \label{M(G,w) = C_0(G, 1/w)* Prop}  For each $\nu \in\M(G,\om)$, $I_\nu \in C_0(G,  \om^{-1} )^*$ and $\|I_\nu\| = \|\nu \|_\om$; moreover, 
\beq \label{M(G,w) = C_0(G, 1/w)* Prop Eqn}  C_0(G,  \om^{-1})^* = \{ I_\nu : \nu \in \M(G, \om)\}. \eeq
We can thus make the identification $\M(G, \om) = C_0(G,  \om^{-1})^* $ through the pairing 
$$\l \nu, \psi \r_\om = \int \psi \, d \nu \qquad (\nu \in \M(G, \om), \psi \in C_0(G,  \om^{-1})).$$
With respect to this identification, the inverse isometric isomorphisms 
$$\M(G, \om) \ra \M_r(G): \nu \mapsto \om \nu  \quad {\rm and } \quad  \M_r(G) \ra \M(G, \om): \mu \mapsto {1 \over \om} \mu$$ 
are weak$^*$-homeomorphisms.  
\ep 

\begin{proof} As noted above, $S: C_0(G,  \om^{-1}) \ra C_0(G): \psi \mapsto {\psi \over \om}$ is an isometric isomorphism, so $S^*: \M_r(G) = C_0(G)^* \ra C_0(G,  \om^{-1})^*$ is also an isometric isomorphism. Let $\nu \in \M(G,\om)$. Then $\om \nu \in  \M_r(G)$ and for $\psi \in C_0(G, \om^{-1})$, ${\psi \over \om} \in C_0(G) \subseteq \L^1(\om \nu)$ (see Remark \ref{M_r=C0* Remark}); therefore by II.7.5, $\psi = (\psi/\om)\om \in \L^1(\nu)$ and 
$$\l I_\nu, \psi \r = \int \psi \, d\nu = \int {\psi \over \om} \, d (\om \nu) = \l \om \nu, S(\psi)\r = \l S^*(\om \nu), \psi \r.$$
Hence, $C_0(G, \om^{-1}) \subseteq \L^1(\nu)$, $I_\nu = S^*(\om \nu) \in C_0(G, \om^{-1})^*$,  and  therefore 
$\|I_\nu\| = \|S^*(\om \nu)\| = \| \om \nu \| = \| \nu\|_\om;$ since $S^*(\mu) = I_{\om^{-1} \mu}$ and $S^*$ maps onto $C_0(G, \om^{-1})^*$, we have (\ref{M(G,w) = C_0(G, 1/w)* Prop Eqn}). Making the identification of $\nu$ and $I_\nu$, $\mu \mapsto {1 \over \om} \mu = S^*(\mu)$ is weak$^*$-continuous, with (weak$^*$-continuous) inverse map $\nu \mapsto \om \nu$.  
\end{proof}

In Lemma \ref{Convolution in M(G,w) Lemma},  $X$ is a locally compact Hausdorff space, $h: X \ra (0, \infty)$ is a continuous function, and $\mu \in \M(X)^+$ is such that $h \mu \in \M_r(X)$. Observe that $\E_\mu \subseteq \fB(X) = \E_{h\mu}$; see Remark \ref{Max regular extension remark}.

\blem  \rm \label{Convolution in M(G,w) Lemma} The function $h$ is locally $\mu_e$-summable and for any set $A \in \E_\mu$, $h(\mu_e)(A) = (h\mu)_e(A)$. 
\elem 

\begin{proof} Let $A \in \E_\mu$.  Take $(C_n)_n$ to be an increasing sequence of compact subsets of $A$ such that $\mu_e(A) = \lim_n \mu(C_n)$ and let $D = \cup_n C_n$. Observe that $D, A\bs D \in \E_\mu$ and 
$\mu_e(D) = \lim \mu_e(C_n) = \lim \mu (C_n) = \mu_e(A)$; hence \beq  \label{Convolution in M(G,w) Lemma Eqn 1}  \mu_e (A \bs D) = 0.  \eeq 
It follows that for any compact subset $C$ of $A \bs D$, $\mu(C)= 0$ and therefore, since $h$ is locally $\mu$-summable and bounded on $C$, $h \mu(C) = 0$.  Hence, 
\beq \label{Convolution in M(G,w) Lemma Eqn 2} \lim (h\mu)_e( A \bs D) = \lim\{ (h \mu)(C): C \in \fC (X), \ C \subseteq A \bs D\} = 0. \eeq 
As noted in Remark \ref{Max regular extension remark}, $h \chi_{A \bs D}$ is locally $\mu_e$-measurable and it follows from (\ref{Convolution in M(G,w) Lemma Eqn 1}) and II.2.7  that 
\beq  \label{Convolution in M(G,w) Lemma Eqn 3}     \int h\chi_{A\bs D} \, d \mu_e = \lim_n \int (h\wedge n) \chi_{A\bs D} \, d \mu_e = 0.   \eeq 
Also, since $h \mu$ is bounded, $ \lim \int h \chi_{C_n} \, d \mu_e = \lim \int h \chi_{C_n} \, d \mu = \sup (h\mu)(C_n) < \infty$
(using II.8.15 Remark 3), and therefore by II.2.7, 
\beq \label{Convolution in M(G,w) Lemma Eqn 4} 
\int h \chi_D \, d\mu_e = \lim \int h \chi_{C_n}  \, d\mu_e = \lim (h\mu)(C_n) = \lim(h\mu)_e(C_n) = (h \mu)_e(D).
\eeq
From (\ref{Convolution in M(G,w) Lemma Eqn 3}) and (\ref{Convolution in M(G,w) Lemma Eqn 4}), $ h \chi_{A\bs D}, h \chi_D \in \L^1(\mu_e)$, whence $h \chi_A \in \L^1(\mu_e)$. Hence, $h$ is locally $\mu_e$ summable. Moreover, (\ref{Convolution in M(G,w) Lemma Eqn 4}), (\ref{Convolution in M(G,w) Lemma Eqn 3})  and (\ref{Convolution in M(G,w) Lemma Eqn 2})  yield $ h(\mu_e)(A) =(h\mu)_e(A)$. 
\end{proof} 

Let $p: G\times G \ra G: (s,t) \mapsto st$. Following III.10.2, we say that $\mu, \nu \in \M(G)$ are convolvable, or that $\mu * \nu$ exists, if $p$ is $\mu \times \nu$-proper in the sense of II.10.3, i.e., if $p^{-1}(A) \in \E_{\mu \times \nu}$ whenever $A \in \fS(G)$. In this case, $\mu * \nu \in \M(G)$, where for $A \in \fS(G)$,  
\beqs  \mu * \nu (A) &=& p_*((\mu \times \nu)_e)(A) = (\mu \times \nu)_e (p^{-1}(A))\\
& = & \lim \{ (\mu \times \nu)(C): C \subseteq p^{-1}(A), \  C \in \fC(G \times G)\};
\eeqs 
see III.10.2, II.10.3, II.10.5, II.10.1.   Equivalently, one can check that $\mu * \nu$ exists if and only if 
$$\sup \{ (|\mu| \times |\nu|)(C) : C \subseteq p^{-1}(D), \  C \in \fC(G \times G) \} < \infty$$ for every compact subset $D$ of $G$.  (In our context, the definition of $\mu \times \nu \in \M(G \times G)$ and its properties are found in \S II.9.) 

\bt  \label{M(G,w) Banach Algebra Thm} \rm  With respect to convolution product, $\M(G, \om)= C_0(G,\om^{-1})^*$ is a Banach algebra, i.e.,   $(\mu, \nu) \mapsto \mu * \nu$ is a well-defined associative operation on $\M(G, \om)$ satisfying $\|\mu* \nu\|_\om \leq \| \mu \|_\om \|\nu\|_\om$. Moreover, for $\mu, \nu \in \M(G, \om)$ and $\psi \in C_0(G, \om^{-1})$,  
\beq \label{Convolution dual pairing description}   \l \mu * \nu, \psi\r_\om = \int \psi(st) \, d(\mu \times \nu)_e(s,t) = \iint \psi(st)\,  d \mu(s) d \nu(t)  =  \iint \psi(st)\,  d \nu(t) d \mu(s).      \eeq  
\et

\begin{proof}  Let $\mu, \nu \in \M(G, \om)$, with $\mu, \nu \geq 0$. Let $D$ be a compact subset of $G$, $C$ a compact subset of $p^{-1}(D)$.  The functions $1_C(x,y)$ and $ g(x,y) = {1 \over \om(x) \om(y)} 1_C(x,y)$ are Borel measurable functions, and are therefore locally $(\sigma \times \rho)$-measurable for any pair of measures $\sigma, \rho \in \M(G)$;  moreover,  since they are non-negative, bounded and vanish off $C$, $1_C, g \in \L^1(\sigma \times \rho)$. Applying the Fubini Theorem (II.9.8) to these functions, and using  II.7.5 twice --- which also applies by II.9.8 --- we obtain
\beqs \mu \times \nu(C) &= &  \int \int 1_C (x,y) \, d\mu(x) d\nu(y) =  \int \int g (x,y) \om(x)  \, d\mu(x) \om(y) \, d\nu(y)\\
& = & \int \int g (x,y) \, d\om \mu(x) \, d\om \nu(y) =  \int_{G \times G} {1 \over \om(x) \om(y)} 1_C (x,y) \, d(\om \mu \times \om \nu)(x,y)\\
& \leq & \int_{G \times G} {1 \over \om(x y)} 1_C (x,y) \, d(\om \mu \times \om \nu)(x,y) \leq \int_{G \times G} M_D  1_C (x,y) \, d(\om \mu \times \om \nu)_e(x,y)
\eeqs
where $ M_D = \sup_{z \in D}  \om(z)^{-1}$, since $C \subseteq p^{-1}(D)$, and we have used II.8.15 Remark 3. Observe that $p^{-1}(D) \in \fB(G \times G) = \E_{\om \mu \times \om \nu}$, since $\om \mu \times \om \nu \in \M_r(G \times G)$ --- see II.9.14 ---  so 
\beqs  \mu \times \nu (C) \leq  \int_{G \times G} M_D  1_{p^{-1}(D)}  \, d(\om \mu \times \om \nu)_e \leq M_D \|\om \mu \times \om \nu \| = M_D \|\om \mu \| \|\om \nu \| = M_D \| \mu \|_\om \| \nu \|_\om.
\eeqs   
Hence, $\mu * \nu$ exists.  We now show  $\mu * \nu \in \M(G, \om)$ and $\| \mu * \nu \|_\om \leq \| \mu\|_\om \| \nu\|_\om$. Let $A \in \fS(G)$.  Since $\om$ is continuous on $G$ and $\mu * \nu \in \M(G)$, $\om$ is locally $\mu * \nu$-summable and $\om(\mu * \nu) \in \M(G)$. Hence, $\om \chi_A \in \L^1(\mu * \nu) = \L^1(p_*(\mu \times \nu)_e)$.  Therefore,  II.10.2 gives $(\om \chi_A) \circ p \in \L^1((\mu \times \nu)_e)$ and 
\beqs \om(\mu * \nu)(A) &=&  \int \om \chi_A \, d (p_*((\mu  \times \nu)_e)) 
 =   \int (\om \chi_A) \circ p  \, d (\mu  \times \nu)_e \\
& =  & \int \om \circ p  \,  \chi_{p^{-1}(A)}  \, d (\mu  \times \nu)_e 
 \leq    \int (\om \times \om)   \chi_{p^{-1}(A)}  \, d (\mu  \times \nu)_e 
\eeqs 
where $(\om \times \om)(s,t) = \om(s) \om(t)$.  By II.9.9 and II.9.3, $(\om \times \om)(\mu \times \nu) = \om \mu \times \om \nu$, which belongs to $\M_r(G \times G)$ by II.9.14. Observe that $\om \times \om$ is locally $(\mu \times \nu)_e$-summable, by Lemma \ref{Convolution in M(G,w) Lemma},  and $p^{-1}(A) \in \E_{\mu \times \nu}$, since $\mu* \nu$ exists. Hence, the above inequality and Lemma \ref{Convolution in M(G,w) Lemma} yield
\beqs \om(\mu * \nu)(A) & \leq & (\om \times \om)(\mu \times \nu)_e(p^{-1}(A)) = ((\om \times \om)(\mu \times \nu))_e(p^{-1}(A)) \\
& = & (\om \mu \times \om \nu)_e(p^{-1}(A)) \leq \| \om \mu \times \om \nu\| =  \| \om \mu\|  \| \om \nu\| = \|  \mu\|_\om  \|  \nu\|_\om.
\eeqs 
Hence, $\om(\mu *\nu)$ is bounded, i.e., $\mu * \nu \in \M(G, \om)$, and 
$\| \mu * \nu \|_\om = \|\om (\mu * \nu)\| \leq \| \mu \|_\om \|\nu\|_\om.$

Assume now that $\mu, \nu $ are any two measures in $\M(G, \om)$. As we have noted, $\sigma \in \M(G, \om)$ exactly when $|\sigma| \in \M(G,\om)$ and $\| \sigma \|_\om = \| | \sigma | \|_\om$, so it follows from III.10.3 and the positive case that $\mu * \nu$ exists and $|\mu * \nu| \leq | \mu | * | \nu |$. Hence, $\om |\mu * \nu | \leq \om | \mu | * | \nu |$, so  $\mu * \nu \in \M(G, \om)$ and 
$$\|\mu * \nu \|_\om = \| \om |\mu * \nu | \| \leq \| \om | \mu | * | \nu | \| = \| | \mu | * | \nu |\|_\om \leq \|| \mu | \|_\om \| | \nu | \|_\om = \| \mu  \|_\om \|  \nu  \|_\om.$$ 
Associativity of convolution in $\M(G, \om)$ is now an immediate consequence of III.10.10.  Since any $\psi \in C_0(G, \om^{-1})$ vanishes off a $\sigma$-compact subset of $G$ and any $\mu, \nu \in \M(G, \om)$ are $\sigma$-bounded  --- since $\om \mu$ and $\om \nu$ are so, and $\om \mu \sim \mu$, $\om \nu \sim \nu$ --- Remark III.10.8 applies to give (\ref{Convolution dual pairing description}). 
\end{proof}

Let $\lambda = \lambda_G$ be a fixed left Haar measure on $G$, $\L^1(G) = \L^1(\lambda)$.  Then $\lambda \in \M(G)$ (\S III.7), so $\om \lambda \in \M(G)$ as well and,  since $\om >0$, $\om \lambda \sim \lambda$, from which it follows that  $g$ is locally $\om \lambda$-measurable and vanishes off a $\sigma$-compact set if and only if $g \om$ is locally $ \lambda$-measurable and vanishes off a $\sigma$-compact set. Hence, if we define $\L^1(G, \om) := \L^1(\om \lambda)$, $g \in \L^1(G, \om)$ exactly when $g\om \in \L^1(G)$, and in this case $ \int g \, d(\om \lambda) = \int g\om \, d \lambda$, by II.7.5. 
Thus, 
$$\L^1(G, \om ) = \{ g : g\om \in \L^1(G) \} \ \ {\rm and} \ \ \|g \|_\om : = \|g \|_{\L^1(\om \lambda)} = \|g\om\|_1$$
defines a Banach space norm on $\L^1(G, \om)$. Moreover, 
$T : \L^1(G, \om ) \ra \L^1(G): g \mapsto g \om$
is an isometric linear isomorphism, with inverse $f \mapsto {1\over \om} f$, so 
$T^*: \L^\infty(G) = \L^1(G)^* \ra \L^1(G, \om)^*$ 
is a weak$^*$-continuous isometric isomorphism given by 
$ \l T^* \phi, g\r_\om = \l \phi, \om g \r = \int (\phi \om) g \, d \lambda.$
Letting 
$$\L^\infty(G, \om^{-1}): = \{ \phi \om: \phi \in \L^\infty(G) \} = \left\{ \psi: {\psi\over \om} \in \L^\infty (G)\right\} \ \ \ {\rm where} \ \ \ \|\psi\|_{\infty, \om^{-1}} := \left\|{\psi\over \om}\right\|_\infty,  $$
 we can hence identify $\L^1(G,\om)^*$ with $\L^\infty(G, \om^{-1})$  via the pairing 
 $\l \psi, g \r_\om = \int \psi g \, d \lambda.$
 Observe that $S = (T^*)^{-1}: \L^\infty(G, \om^{-1}) \ra \L^\infty(G): \psi \mapsto {\psi \over \om}$
 is a weak$^*$-homeomorphic isometric isomorphism. 
(We note that $\L^\infty(G, \om^{-1})$ is not usually the same space as $\L^\infty(\om \lambda)$ ($= \L^\infty(\lambda)$ because $\om \lambda \sim \lambda$), which can also be identified with $\L^1(\om \lambda)^* = \L^1(G, \om)^*$ in the usual way by II.7.11.) Note that because $T^{-1}$ maps $C_{00}(G)$ onto itself, $C_{00}(G)$ is dense in $\L^1(G,\om)$. 

Let $g \in \L^1(G,\om) = \L^1(\om \lambda)$, $A \in \fS(G)$. Then $\om g \in \L^1(\lambda)$ and ${1\over \om}$ is bounded on $A$, so $\chi_A g = ({1\over \om} \chi_A) \om g \in \L^1(\lambda)$; hence,   $g \lambda \in \M(G)$ is well-defined (II.7.2). Also, $\om (g \lambda) = (\om g) \lambda \in \M(G)$ by II.7.5 and, by II.7.9/III.11.3, $\|f \|_1 = \|f \lambda\|$ for $f \in \L^1(G)$ and  
$$\M_a(G) = \{ \mu \in \M_r(G): \mu \ll \lambda \} = \{ f \lambda : f \in \L^1(G)\} = \{ (\om g) \lambda: g \in \L^1(G, \om)\}.  $$
Since $\om \nu \sim \nu$ for any $\nu \in \M(G)$, it readily follows that  $g \mapsto g\lambda: \L^1(G,\om) \ra \M_a(G, \om)$ is a surjective linear isometry, where 
$\M_a(G, \om) : = \{ \nu \in \M(G, \om) : \nu \ll \lambda \}. $
We can thus identify $\L^1(G, \om)$ with $\M_a(G, \om)$ via $g \mapsto g \lambda$. 

\bp \label{L1(G,w) Ideal Prop} \rm The Banach space $\L^1(G, \om) = \M_a(G, \om)$  is a closed ideal in $\M(G, \om)$ and has a positive contractive approximate identity. Moreover, if $g\in \L^1(G, \om)$ and $\nu \in \M(G,\om)$, then $\nu * g, g* \nu \in \L^1(G, \om)$ are given by the formulas, which hold for locally $\lambda$-almost all $t \in G$,
\beq   \label{L1(G,w) Ideal Prop Eqn1} \nu * g(t) = \int g(s^{-1}t) \, d \nu(s) \quad {\rm and } \quad g * \nu(t) = \int \Delta(s^{-1}) g(ts^{-1}) \, d \nu (s);
 \eeq
thus, $\L^1(G, \om)$ is a Banach algebra with respect to the convolution product 
\beq   \label{L1(G,w) Ideal Prop Eqn2} f*g(t) = \int f(s) g(s^{-1}t) \, d \lambda(s). 
 \eeq
 \ep 

\begin{proof}  We have already noted that $g$ is locally $\lambda$-summable and vanishes off a $\sigma$-compact set, and $(g \lambda)* \nu$, $\nu * (g \lambda)$ exist in $\M(G, \om)$ by Theorem \ref{M(G,w) Banach Algebra Thm}. Letting $h(t)$ and $k(t)$ be defined by the respective integral formulas on the left and right of (\ref{L1(G,w) Ideal Prop Eqn1}), $\nu * (g \lambda) = h \lambda$ and $(g \lambda) * \nu = k \lambda$ by III.11.5. Thus, $h \lambda, k \lambda \in \M_a(G, \om) = \{ f \lambda : f \in \L^1(G, \om)\}$, so the uniqueness part of the Radon--Nikodym Theorem --- see Remark 1 of II.7.8 --- implies that $h, k \in \L^1(G, \om)$. The formula (\ref{L1(G,w) Ideal Prop Eqn2}) now follows quickly (or directly from III.11.6). Let $\cal I$ be the  neighbourhood system at $e_G$ and for each $\alpha \in {\cal I}$, let $f_\alpha \in C_{00}(G)$ be chosen with $f_\alpha \geq 0$, $\|f_\alpha\|_1 = 1$ and   support contained in $\alpha$.  Then $(f_\alpha)_\alpha$ is a bounded approximate identity for $\L^1(G)$.  Letting $ e_\alpha = \om^{-1} f_\alpha$, $\|e_\alpha\|_{\om} =1$ and $\|e_\alpha\|_1 \ra 1$, from which it easily follows that $(e_\alpha)_\alpha$ is  also a bounded approximate identity for $\L^1(G)$; the proof of Lemma 2.1 in \cite{Gha1} now shows that $(e_\alpha)_\alpha$ is a contractive approximate identity for $\L^1(G, \om)$. 
\end{proof} 

\br \rm  Every Borel measurable function is locally $\lambda$-measurable and  every $f \in L^1(G, \om)$ --- where $L^1(G,\om)$ is defined in the usual sense (as in the introduction) --- vanishes off a $\sigma$-compact set. It follows that the Banach algebra $\L^1(G, \om)$, as we have defined it, exactly coincides with the usual definition of the Beurling group algebra $L^1(G,\om)$, which, as noted in the introduction, is always valid.  Going forward, we can therefore use any known result about $L^1(G, \om) = \L^1(G, \om)$ that was proved independently of $M(G, \om)$. 
 \er

\section{The dual Banach algebra $\M(G,\om)$  and the embedding map } 

The \it support \rm of  $\mu$ in   $\M(G)$ is the set ${\rm s}(\mu) = G \bs \bigcup \{ U \in \fS(G): U \ {\rm is \ open \ and \ } |\mu|(U)= 0 \} $(II.8.9). Let 
$\M_{cr}(G) = \{ \mu \in \M(G): {\rm s}(\mu) {\rm \ is \ compact}\}.$

\br \rm 1. Observe that ${\rm s}(\mu) =  {\rm s}(\mu_e) = G \bs \bigcup \{ V \in \E_\mu :  V \ {\rm is \ open \ and \ } |\mu_e|(V)= 0 \} $. 

\smallskip 

\noindent  2.  Since $\om$ and ${1\over \om}$ are bounded on any set $A$ in $\fS(G)$,    ${\rm s}(\mu) = {\rm s}(\om \mu) = {\rm s}\left({1\over \om}\mu\right)$ for any $\mu \in \M(G)$. 

\smallskip 

\noindent  3. By III.10.16, $\M_{cr}(G)$ is a dense subalgebra of $\M_r(G)$.   From 2 above, the inverse linear isometries $\nu \mapsto \om \nu$ and $\mu \mapsto {1\over \om} \mu$ between $\M(G, \om)$ and $\M_r(G)$ map $\M_{cr}(G)$ onto itself, so $\M_{cr}(G)$ is also a  dense subalgebra of $\M(G, \om)$.  

\er 

A measure $\sigma$ on a $\delta$-ring $\fS$ is \it concentrated \rm on a set $F$  if for each $A \in \fS$, $A \cap F, A \bs F \in \fS$ and $\sigma (A) = \sigma (A \cap F)$ or, equivalently, $\sigma(A \bs F) = 0$.  For $\mu \in \M(G)$ and  a Borel set $F$, $A \cap F, A \bs F \in \fS(G) $ (respectively, $A \cap F, A \bs F \in {\cal E}_\mu$) is automatic for any $A \in \fS(G)$ ($A \in {\cal E}_\mu$), and it is clear from (\ref{Max regular extension eqn})  that $\mu$ is concentrated on $F$ if and only if $\mu_e$ is concentrated on $F$.  A function $\psi \in LUC(G, \om^{-1})$  may fail to vanish off a $\sigma$-compact set and therefore, as noted in Remark \ref{M_r=C0* Remark}, in this theory we cannot integrate $\psi$  with respect to any $\mu $ in $\M(G)$. Lemma \ref{Concentrated Measures Lemma} allows us to move past this issue.  

\blem  \label{Concentrated Measures Lemma} \rm (a) Every  $\mu$ in  $\M(G)$ is concentrated on its support, ${\rm s}(\mu)$. \\
  (b) Let $\mu \in \M_r(G)$.  Then $\mu$ (and therefore $\mu_e$) is concentrated on a $\sigma$-compact subset $F$ of $G$ and, for any such $F$ and any  Borel measurable function $f \in \L^1(\mu_e)$,   $f \chi_F \in \L^1(\mu)  \cap \L^1(\mu_e)$ and $$\int f \, d \mu_e =  \int f \chi_F \, d \mu_e = \int f \chi_F \, d \mu.$$
(c) Any $\nu \in \M(G, \om)$ is concentrated on a $\sigma$-compact set.   
 \elem

\begin{proof}   (a) Let  $A \in \fS(G)$.  Any  compact subset of   $A \bs {\rm s}(\mu)$ is covered by   the collection  of open sets  $U \in \fS(G)$ with  $|\mu|(U)= 0$, and is therefore    $|\mu|$-null; by regularity of $\mu$ (II.8.2(II)),  $|\mu|(A \bs {\rm s}(\mu)) = 0$.

\noindent (b)  Take $(C_n)_n$ to be an increasing sequence of compact subsets of ${\rm s}(\mu)$ such that $|\mu|(C_n) > \|\mu \|-1/n$ and let $F = \bigcup C_n$, where we have used (b).  Then $\mu$ is concentrated on $F$ because for $A \in \fS(G)$, 
$$|\mu|(A \bs F) = |\mu| ((A \bs F)\cap {\rm s}(\mu)) \leq |\mu|_e ({\rm s}(\mu) \bs F)  = |\mu|_e({\rm s}(\mu)) - |\mu_e|(F) = \| \mu \| - \lim |\mu|(C_n) = 0.$$
Suppose  $\mu \geq 0$, $F$ is  any $\sigma$-compact set on which $\mu $ is concentrated, and $f \in \L^1(\mu_e)$ is a non-negative Borel-measurable function.  It is then clear (from II.2.2 and II.2.5) that $f \chi_F \in \L^1(\mu_e)$ and $ \int f  \, d \mu_e = \int f \chi_F \, d \mu_e$. Also, $f \chi_F $ is locally $\mu$-measurable (II.8.2), vanishes off the  $\sigma$-compact set $F$ and, taking any sequence of non-negative $\fS(G)$-simple functions such that $h_n \uparrow f \chi_F$,  II.2.2 gives 
$$\int f \chi_f \, d\mu_e = \lim \int h_n \, d \mu_e = \lim \int h_n \, d \mu = \int f \chi_F \, d \mu.$$ 
(c) Since $\om \nu \in \M_r(G)$ and $\nu \sim \om \nu$, this follows from (b). 
 \end{proof}

Since $\L^1(G, \om)$ is a closed ideal in $\M(G, \om)$,  through the operations 
$$\l \psi \c \nu, g \r = \l \psi, \nu * g \r \ \ {\rm and} \ \ \l \nu \c \psi, g \r = \l \psi, g* \nu  \r  \quad (\psi \in \L^\infty(G, \om^{-1}), \nu \in \M(G, \om), g\in \L^1(G, \om)), $$
$\L^\infty(G, \om^{-1}) = \L^1(G, \om)^*$ is a dual $\M(G,\om)$-module.
Observe that for $\psi \in \L^\infty(G, \om^{-1})$  and $s \in G$, 
$$ \psi \c \delta_s (t) = \psi \c s (t) := \psi(st) \ \ {\rm and } \ \  \delta_s \c \psi (t) = s \c \psi (t) := \psi(ts) \quad (t \in G).$$
Recall that $\psi$ belongs to  $LUC(G, \om^{-1})$ [$RUC(G, \om^{-1})$] when $ {\psi \over \om}$ belongs to $LUC(G)$ [$RUC(G)$]. For $LUC(G, \om^{-1})$, the following is \cite[Proposition 1.3]{Gro} and \cite[Propositions 7.15 and 7.17]{Dal-Lau}, (where no restrictions are needed on the weight $\om$); symmetric arguments establish the $RUC(G, \om^{-1})$ case.

\blem \label{Gronbaek Prop} \rm The following statements are equivalent: \bi \item[(a)]   $\psi \in LUC(G, \om^{-1})$ [$RUC(G, \om^{-1})$];
\item[(b)] $\psi \in \ell^\infty(G, \om^{-1})$ and the map $G \ra (\ell^\infty(G, \om^{-1}), \| \c \|_{\infty, \om^{-1}}): s \mapsto \psi \c s \ [s \c \psi]$ is continuous;
 \item[(c)] $\psi \in \L^\infty(G, \om^{-1})$ and  the map $G \ra (\L^\infty(G, \om^{-1}), \| \c \|_{\infty, \om^{-1}}): s \mapsto \psi \c s \ [s \c \psi]$ is continuous; 
 \item[(d)] $\psi \in \L^\infty(G, \om^{-1}) \c \L^1(G, \om)$ [$\psi \in \L^1(G, \om) \c \L^\infty(G, \om^{-1})$]. 
\ei
\elem

\br  \rm 1. Observe that  condition (b)  implies  $\psi $ is continuous on $G$, whence $\psi \in \L^\infty(G, \om^{-1})$. 

\smallskip 

\noindent 2. In the proof of \cite[Proposition 7.15]{Dal-Lau}, the authors establish continuity of a function $\psi$ satisfying (c) via Ascoli's theorem.  An alternative approach is to  establish (i) and (ii) as follows:\\
 (i) If $\phi \in \L^\infty(G, \om^{-1})$ and $g \in \L^1(G, \om)$, then $\phi \c g$ can be identified with the continuous function  \beq \label{First Module Action LUC Eqn} (\phi \c g)(t) = \l \phi, g* \delta_t\r \qquad {\rm for \ every } \ t \in G. \eeq 
 [Note that $H \in \ell^\infty(G, \om^{-1})$ where $H(t):= \l \phi, g* \delta_t\r$ and, since    $t \mapsto g * \delta_t: G \ra (\L^1(G, \om), \|\c \|_\om)$ is  continuous --- e.g., see \cite[Lemma 3.1.5]{Zad0}, which holds for any weight $\om$ --- $H$ is continuous on $G$ (and satisfies Lemma \ref{Gronbaek Prop}(c)); in a standard way, one can  check that for $f \in \L^1(G, \om)$, $\l \phi \c g, f \r = \l H, f\r$.] \\
  (ii) If $\psi$ satisfies (c) and $(e_i)$ is a bounded approximate identity for $\L^1(G, \om)$, then $\| \psi\c e_i  - \psi\|_{\infty, \om^{-1}} \ra 0$; since $CB(G, \om^{-1})$ is closed in $\L^\infty(G, \om^{-1})$, $\psi \in CB(G, \om^{-1})$. 
\er 

\bp   \label{LUC(G,1/w) a M(G,w) module Prop}  \rm The spaces $LUC(G, \om^{-1})$ and $RUC(G, \om^{-1})$ are $\M(G, \om)$-submodules of $\L^\infty(G, \om^{-1})$. Moreover, for $\nu \in \M(G, \om)$, $\psi \in LUC(G, \om^{-1})$ [$\psi \in RUC(G, \om^{-1})$] and for every $s \in G$, 
$$(\nu \c \psi)(s) = \int (\psi \c s) \chi_{F_s} \, d \nu = \int {\psi \c s \over \om}  \, \chi_{F_s} \, d (\om \nu) = \int {\psi \c s \over \om}   \, d (\om \nu)_e$$
$$ \left[(\psi\c \nu )(s) = \int (s \c \psi ) \chi_{F_s} \, d \nu = \int {s \c \psi  \over \om}  \, \chi_{F_s} \, d (\om \nu) = \int {s \c \psi  \over \om}   \, d (\om \nu)_e\right],$$
where $F_s$ is any $\sigma$-compact set on which $\nu$ is concentrated; $F_s$ can be chosen to vary with $s\in G$. 
\ep

\begin{proof}   Letting $\nu \in \M(G, \om)$,  $\psi \in LUC(G, \om^{-1})$, it is clear from Lemma \ref{Gronbaek Prop} (d) that $\psi \c \nu, \nu \c \psi \in LUC(G, \om^{-1})$.  Since ${\psi \c s \over \om} \in LUC(G)$ and $\om \nu \in \M_r(G)$, 
$$H(s) = H_{\nu, \psi}(s) := \int {\psi \c s \over \om}   \, d (\om \nu)_e = \int {\psi \c s \over \om}  \, \chi_{F_s} \, d (\om \nu)$$ is well-defined, where we have used Lemma \ref{Concentrated Measures Lemma}.  The function $(\psi \c s) \chi_{F_s} \in \ell^\infty(G, \om^{-1})$ is Borel measurable --- and therefore locally $\nu$-measurable --- and vanishes off the $\sigma$-compact set $F_s$, so 
$ {\psi \c s \over \om}  \chi_{F_s} \in \L^1(\om \nu)$. Therefore, by II.7.5, $(\psi \c s ) \chi_{F_s} \in \L^1(\nu)$ and 
$$\int (\psi \c s) \chi_{F_s} \, d \nu  =  \int {\psi \c s\over \om} \,  \chi_{F_s}  \, \om \, d \nu  = \int {\psi \c s \over \om}  \, \chi_{F_s} \, d (\om \nu) = H(s).$$
Since $\ds |H(s)| \leq \left\|{\psi\c s\over \om}\right\|_\infty\|\om \nu\| \leq \om(s) \|\psi\|_{\infty, \om^{-1}} \|\nu \|_\om$, $H = H_{\nu, \psi}\in \ell^\infty(G, \om^{-1})$ with $\| H_{\nu, \psi}\|_{\infty, \om^{-1}} \leq \|\psi\|_{\infty, \om^{-1}} \| \nu\|_\om$. Hence, if $s_i \ra s$ in $G$, 
$$\|(H_{\nu, \psi}) \c s_i - (H_{\nu, \psi})\c s\|_{\infty, \om^{-1}} = \|H_{\nu, \psi \c s_i - \psi \c s} \|_{\infty, \om^{-1}} \leq \|\psi  \c s_i  - \psi \c s\|_{\infty, \om^{-1}} \| \nu \|_\om \ra 0; $$
by Lemma \ref{Gronbaek Prop}, $H_{\nu, \psi} \in LUC(G, \om^{-1})$.
To show that $H_{\nu, \psi} = \nu \c \psi$, we can assume $\nu \geq 0$,  $\psi\geq 0$ and take $F = F_s$ for each $s \in G$.  Let $g\geq 0$ be a function in the dense subspace $C_{00}(G)$ of  $\L^1(G, \om)$.  Since the maps $(s,t) \mapsto \psi(t) \Delta(s^{-1}) g(ts^{-1}) \chi_F(s), \ \psi(ts)g(t) \chi_F(s)$ are Borel measurable --- hence locally $(\nu \times \lambda)$-measurable --- and vanish off a $\sigma$-compact subset of $G \times G$, our applications of the Fubini Theorem (II.9.8) are valid in the following calculation. Using (\ref{L1(G,w) Ideal Prop Eqn1}): 
\beqs 
\l \nu \c \psi, g \r & = & \l \psi, g * \nu \r = \int \psi(t) \int  \Delta(s^{-1}) g(ts^{-1})\, d \nu(s) \, d \lambda (t)  \\ 
&= &  \int \int \psi(t)  \Delta(s^{-1}) g(ts^{-1}) \chi_F(s) \, d \nu(s) \, d \lambda (t) 
 =   \int  \int \psi(t)    \Delta(s^{-1}) g(ts^{-1}) \chi_F(s) \, d \lambda(t) \, d \nu (s) \\
 & = &   \int  \int \psi(ts)  g(t) \chi_F(s) \, d \lambda(t) \, d \nu (s) 
  =   \int  \int \psi\c t(s)  \chi_F(s) \, d \nu(s) \, g(t)  \, d \lambda(t)  =   \l H_{\nu,\lambda}, g \r; 
\eeqs
since both functions are continuous, $\nu \c \lambda = H_{\nu, \lambda}$. 
\end{proof}

\bc \label{C_0(G,1/w) a M(G,w) module Corollary} \rm  The space $C_0(G, \om^{-1})$ is a  $\M(G, \om)$-submodule of $\L^\infty(G, \om^{-1})$, and  for $\nu \in \M(G, \om)$,  $\psi \in C_0(G, \om^{-1})$ and $s \in G$, 
\beq \label{C_0(G,1/w) a M(G,w) module Corollary Eqn} \nu \c \psi(s) = \int \psi \c s  \, d \nu = \l \nu, \psi \c s \r_\om \ \ {\rm and } \ \  \psi\c \nu (s) = \int s \c \psi  \, d \nu = \l \nu,  s \c \psi\r_\om. \eeq 
\ec

\begin{proof}  Let $\psi \in C_0(G, \om^{-1})$ and let $F$ be a $\sigma$-compact set on which $\nu$ is concentrated. Taking $A_s$ to be a $\sigma$-compact set off of which $\psi \c s$ and $s \c \psi$ vanish, and putting $F_s = F \cup A_s$, Proposition \ref{LUC(G,1/w) a M(G,w) module Prop} gives $\nu \c \psi, \psi \c \nu \in (LUC\cap RUC)(G, \om^{-1})$ and 
$$\nu \c \psi(s) = \int (\psi \c s) \chi_{F_s} \, d \nu = \int \psi \c s  \, d \nu  \ \ {\rm and } \ \  \psi\c \nu (s) = \int (s \c \psi ) \chi_{F_s} \, d \nu  = \int s \c \psi  \, d \nu. $$
Observe that $\nu \c \psi$ is supported on ${\rm s}(\psi) {\rm s}(\nu)^{-1}$, which is compact when $\nu$ belongs to the dense subspace $\M_{cr}(G)$ of $\M(G, \om)$ and $\psi $ belongs to the dense subspace $C_{00}(G)$ of $C_0(G, \om^{-1})$. It follows that $C_0(G, \om^{-1})$ is a left (and similarly, right) $\M(G, \om)$-submodule of $\L^\infty(G, \om^{-1})$. 
\end{proof}

It follows that $\M(G,\om) = C_0(G, \om^{-1})^*$ is a dual $\M(G, \om)$-module with respect to the operations
$$\l \mu \c_r \nu, \psi \r_\om = \l \mu,  \nu \c \psi \r_\om \quad   {\rm and} \quad  \l \mu \c_{l} \nu, \psi \r_\om = \l  \nu, \psi \c \mu \r_\om \qquad (\mu, \nu \in \M(G, \om), \psi \in C_0(G, \om^{-1})).$$ However, from (\ref{Convolution dual pairing description}) and (\ref{C_0(G,1/w) a M(G,w) module Corollary Eqn}), \beq \label{Convolution as dual module Eqn} \mu\c_r \nu = \mu * \nu = \mu \c_l \nu,\eeq so  $(\mu, \nu) \mapsto \mu  * \nu$ is separately weak$^*$-continuous on $\M(G, \om)$. Hence: 

\bc \rm The Beurling measure algebra $\M(G, \om)$ is a dual Banach algebra. 
\ec 

Let $A$ be a Banach algebra. Recall that a closed submodule ${\cal S}(A^*)$ of the dual $A$-bimodule $A^*$ is left [right] introverted if for each $\mu \in {\cal S}(A^*)^*$ and $\phi \in {\cal S}(A^*)$, $ \mu \sq \phi \in{\cal S}(A^*) $ [$\phi \di \mu \in {\cal S}(A^*)$] where $\mu \sq \phi, \phi \di \mu \in A^*$ are defined by 
$$  \l \mu \sq \phi, a \r_{A^*-A} = \l \mu, \phi \c a \r_{{\cal S}^*-{\cal S}} \quad {\rm and} \quad  \l \phi \di \mu, a \r_{A^*-A} = \l \mu,  a \c \phi \r_{{\cal S}^*-{\cal S}};$$
in this case, ${\cal S}(A^*)^*$ is a  Banach algebra with respect  to its left [right] Arens product 
$$\l \mu \sq \nu, \phi \r = \l \mu, \nu \sq \phi \r \quad [\l \mu \di \nu, \phi \r = \l \nu, \phi \di \mu \r] \qquad (\mu, \nu \in {\cal S}(A^*)^*, \ \phi \in {\cal S}(A^*)).$$
The map $\eta_{\cal S}: A \ra {\cal S}(A^*)^*$ defined by $\l \eta_{\cal S}(a), \phi\r = \l \phi, a\r $ is a bounded homomorphism with weak$^*$-dense range  and,  when $A$ is left introverted,  $\eta_{\cal S}$ maps into the topological centre of $({\cal S}(A^*)^*, \sq)$, $Z_t({\cal S}(A^*)^*) = \{ \mu \in {\cal S}(A^*)^*: \nu \mapsto \mu \sq \nu \ {\rm is \ wk}^*-{\rm wk}^* \ {\rm continuous \ on \ } {\cal S}(A^*)^*\}.$
 For  this see, e.g., \cite{Dal-Lau}.  

\bigskip 

\bp \rm The subspace $C_0(G, \om^{-1})$ of $\L^\infty(G, \om^{-1}) = \L^1(G, \om)^*$ is left and right introverted and $\mu * \nu = \mu \sq \nu = \mu \di \nu$  for $\mu, \nu \in \M(G, \om)= C_0(G, \om^{-1})^*$.  \ep 

\begin{proof}  By Corollary \ref{C_0(G,1/w) a M(G,w) module Corollary}, $C_0(G, \om^{-1})$ is a $\L^1(G, \om)$-submodule of $\L^\infty(G, \om^{-1})$. Let $\mu, \nu \in \M(G, \om)$, $\psi \in C_0(G, \om^{-1})$.  For $g \in \L^1(G, \om)$, equation  (\ref{Convolution as dual module Eqn}) gives
\beqs  \l \nu \sq \psi, g\r_{\L^\infty - \L^1} = \l \nu, \psi \c g \r_\om =  \l g* \nu, \psi \r_\om = \l g, \nu \c \psi \r_\om = \l \nu \c \psi, g \r_{\L^\infty - \L^1}.
\eeqs 
 Hence, $C_0(G, \om)$ is left introverted and   
$\l \mu \sq \nu, \psi \r = \l \mu, \nu \sq \psi \r = \l \mu, \nu \c \psi\r = \l \mu * \nu, \psi\r$, where we have again used (\ref{Convolution as dual module Eqn}).
 Similarly, $C_0(G, \om^{-1})$ is right introverted and $\mu * \nu = \mu \di \nu$. 
\end{proof}    

Let ${\cal S}(\om^{-1})$ be a left introverted subspace of $\L^\infty(G, \om^{-1})$ such that $C_0(G, \om^{-1}) \preceq {\cal S}(\om^{-1}) \preceq LUC(G,\om^{-1})$ and define  
\beq \label{Def of embedding map theta} \Theta: \M(G, \om) \ra {\cal S}(\om^{-1})^* \ \ {\rm by} \ \ \l \Theta(\nu), \psi\r_{{\cal S}^* - {\cal S}} = (\nu \c \psi)(e_G) = \int \psi \, \chi_{F_\nu} \, d \nu \eeq  where $\nu \in \M(G, \om)$, $\psi \in {\cal S}(\om^{-1})$  and 
 $F_\nu$ is any $\sigma$-compact set on which $\nu$ is concentrated.   By Proposition \ref{LUC(G,1/w) a M(G,w) module Prop}, $\Theta$ is well-defined and 
$|\l \Theta(\nu), \psi\r| \leq \|\nu \c \psi \|_{\infty, \om^{-1}} \leq \| \nu \|_\om \|\psi\|_{\infty, \om^{-1}}$,  so $\|\Theta (\nu)\| \leq \| \nu \|_\om$; by equation (\ref{C_0(G,1/w) a M(G,w) module Corollary Eqn}), $\Theta(\nu)\Large{|}_{C_0(G, \om^{-1})} = \nu$, so $\| \Theta(\nu) \| = \|\nu\|_\om$.  Thus, $\Theta$ is a linear isometry.  

Let $so_l$
and $so_r$ denote  the left and right strict
topologies on $\M(G, \om)$ taken with respect to the ideal $\L^1(G, \om)$, i.e.,   the locally convex topologies  respectively generated by the semi-norms
$p_g(\nu) = \|g * \nu \|$ and $q_g(\nu)= \| \nu * g \|$ for $g \in \L^1(G, \om), \nu \in
\M(G, \om)$. Since $\L^1(G, \om)$ has a contractive approximate identity,  (the unit ball of) $\L^1(G,\om)$ is $so_l/so_r$-dense in (the unit ball of) $\M(G, \om)$. 
Observe that when ${\cal S}(\om^{-1})\preceq LUC(G, \om^{-1})$ is a $\L^1(G, \om)$-submodule of $\L^\infty(G, \om^{-1})$, by Lemma \ref{Gronbaek Prop}(d)  
and the  Cohen factorization theorem \cite[Theorem 11.10]{Bon-Dun}, ${\cal S}(\om^{-1}) = {\cal S}(\om^{-1}) \c \L^1(G, \om)$.  Also note that $LUC(G, \om^{-1})$ is always left introverted in $\L^\infty(G, \om^{-1})$ by Lemma \ref{Gronbaek Prop} and \cite[Proposition 5.9]{Dal-Lau}.
 In the non-weighted case and when $\om \geq 1$, the final statement in Proposition \ref{Strict-wk* conts embedding Prop}, which simplifies Arens product calculations, is \cite[Lemma 3]{Lau} and \cite[Proposition 7.21]{Dal-Lau}, respectively.   
 
\bp  \rm  \label{Strict-wk* conts embedding Prop}  Suppose that ${\cal S}(\om^{-1})$ is a left [right] introverted subspace of $\L^\infty(G, \om^{-1}) = \L^1(G, \om)^*$ and $C_0(G, \om^{-1}) \preceq {\cal S}(\om^{-1}) \preceq LUC(G, \om^{-1}) \  [RUC(G,\om^{-1})]$. Then $\Theta: \M(G, \om) \hookrightarrow {\cal S}(\om^{-1})^*$ is a $so_l$-weak$^*$ [$so_r$-weak$^*$] continuous isometric homomorphic embedding into $Z_t({\cal S}(\om^{-1})^*)$ that extends $\eta_{\cal S}: \L^1(G, \om) \ra {\cal S}(\om^{-1})^*$. Moreover,  $(n \sq \psi)(s) = \l n, \psi \c s\r$ for any  $n  \in {\cal S}(\om^{-1})^*$, $\psi \in {\cal S}(\om^{-1})$ and $s \in G$; hence,  ${\cal S}(\om^{-1})$ is introverted as a subspace of $\ell^\infty(G, \om^{-1}) = \ell^1(G, \om)^*$,   the  Arens product on ${\cal S}(\om^{-1})^*$ agrees under either interpretation, and $\Theta$ also extends $\eta_{\cal S}: \ell^1(G, \om) \hookrightarrow {\cal S}(\om^{-1})^*$.   
 \ep 
 
 \begin{proof}    If $g \in \L^1(G, \om) = \L^1(\om \lambda)$, $g$ vanishes off a $\sigma$-compact set $F_g$, and therefore $g = g\lambda \in \M(G, \om)$ is concentrated on $F_g$; hence, for $\psi \in {\cal S}(\om^{-1})$, 
 $$\l \Theta(g), \psi \r_{{\cal S}^*-{\cal S}} = \int \psi \, \chi_{F_g} \, d(g\lambda) =  \int \psi g  \, d\lambda = \l \psi, g\r_{\L^\infty - \L^1} = \l \eta_{\cal S}(g), \psi\r_{{\cal S}^*-{\cal S}}.$$
 For $f \in \L^1(G, \om)$, $\nu \in \M(G, \om)$ and $\psi \in {\cal S}(\om^{-1})$, 
 \beqs  \l \Theta(f) \sq \Theta(\nu), \psi\r_{{\cal S}^*-{\cal S}} & = & \l \eta_{\cal S} (f), \Theta(\nu) \sq \psi\r_{{\cal S}^*-{\cal S}} = \l \Theta(\nu) \sq \psi, f\r_{\L^\infty - \L^1} 
  =  \l \Theta(\nu), \psi  \c f \r_{{\cal S}^*-{\cal S}}  \\  & = &  \nu \c (\psi \c f)(e_G)   =  (\nu \c \psi) \c f(e_G) 
 =  \l \nu \c \psi, f* \delta_{e_G} \r_{\L^\infty - \L^1} \\
 &  =  &  \l \psi, f* \nu \r_{\L^\infty - \L^1} = \l \eta_{\cal S} (f * \nu), \psi \r_{{\cal S}^*-{\cal S}} = \l \Theta(f * \nu), \psi \r_{{\cal S}^*-{\cal S}}, 
 \eeqs
 where we have used (\ref{First Module Action LUC Eqn}). Suppose  that $\nu_i \ra \nu $ $so_l$.  Writing $\psi \in {\cal S}(\om^{-1})$ as $\psi = \phi \c g$ for some $\phi \in {\cal S}(\om^{-1})$ and $g \in \L^1(G, \om)$, 
 \beqs \l \Theta(\nu_i) - \Theta(\nu), \psi\r_{{\cal S}^*-{\cal S}} 
  =   \l \Theta (g)  \sq  \Theta (\nu_i -  \nu),  \phi\r_{{\cal S}^*-{\cal S}}  =  \l \Theta (g * (\nu_i -  \nu)),  \phi\r_{{\cal S}^*-{\cal S}} \ra 0.
 \eeqs
 Hence, $\Theta$ is $so_l$-weak$^*$ continuous. Let  $\mu, \nu \in \M(G, \om)$ and let $(h_i)$ be a net in $\L^1(G, \om)$ such that $so_l-\lim h_i = \mu$. Then $so_l-\lim h_i * \nu   = \mu * \nu$, so 
 $$\Theta(\mu) \sq \Theta(\nu) = {\rm wk}^*-\lim \Theta(h_i) \sq \Theta(\nu) = {\rm wk}^*-\lim \Theta(h_i * \nu) = \Theta(\mu * \nu). $$
 Identify  the Banach algebra $\M(G, \om)$ with its copy $\Theta(\M(G, \om))$ in ${\cal S}(\om^{-1})^*$. Since ${\cal S}(\om^{-1}) = {\cal S}(\om^{-1}) \c \L^1(G, \om)$ is a right $\M(G, \om)$-module, ${\cal S}(\om^{-1})^*$ is a left dual $\M(G, \om)$-module, and the proof of \cite[Lemma 1.4]{Gha-Zad} shows that $\mu \sq n = \mu \c n$ for $\mu \in \M(G, \om)$ and $n \in {\cal S}(\om^{-1})^*$; hence, $\Theta$ maps into $Z_t({\cal S}(\om^{-1})^*)$. For $n  \in {\cal S}(\om^{-1})^*$, $\psi \in {\cal S}(\om^{-1})$ and $s \in G$, 
 $  (n \sq \psi)(s) = \l \delta_s, n \sq \psi\r = \l \delta_s \sq n, \psi\r= \l \delta_s \c n, \psi \r = \l n, \psi \c \delta_s\r= \l n, \psi \c s\r. $  
 The final line is now easily verified. \end{proof} 

 For a Banach algebra $A$, the space $WAP(A^*)$ of weakly almost periodic functionals on $A$ is a left and right introverted subspace of $A^*$ such that for every $m,n \in WAP(A^*)^*$, $m \sq n = m \di n$  \cite[Proposition 3.11]{Dal-Lau}.  Thus, $WAP(A^*)^*$ is a dual Banach algebra. Moreover,  $WAP(A^*)^*$  satisfies the following universal property \cite[Theorem 4.10]{Run}.
 
 \bt \label{Runde Thm} (Runde) \rm
  If $\fB$ is a dual Banach algebra and  $\vp : A \ra \fB$ is a continuous algebra homomorphism, then there is a unique weak$^*$-weak$^*$ continuous algebra homomorphism $\vp_{WAP}: WAP(A^*)^* \ra \fB$ such that $\vp_{WAP} \circ \eta_{WAP} = \vp$. 
\et 

Taking $A_\om = \L^1(G, \om)$, it follows that the embedding ${\rm id}: \L^1(G, \om) \hookrightarrow \M(G, \om)$ determines a unique weak$^*$-weak$^*$ continuous homomorphism $P: WAP(A_\om^*)^* \ra \M(G, \om)$ such that $P \circ \eta_{WAP} = {\rm id}$.  Letting $P_*: C_0(G, \om^{-1}) \ra  WAP(A_\om^*) \preceq \L^\infty(G, \om^{-1})$ denote the predual mapping of $P$, 
$\l P_* \psi, g\r_{\L^\infty-\L^1} = \l P \circ \eta_{WAP} (g), \psi \r = \l \psi, g\r_{\L^\infty - \L^1}$ for $\psi \in C_0(G, \om^{-1})$, $g \in \L^1(G, \om)$. Hence, $C_0(G, \om^{-1} ) \preceq WAP(A_\om^*)$. Moreover, by \cite[Proposition 3.12]{Dal-Lau} and Lemma \ref{Gronbaek Prop}, $WAP(A_\om^*) \preceq (LUC \cap RUC)(G, \om^{-1})$.  Hence, we have the following immediate corollary to Proposition \ref{Strict-wk* conts embedding Prop}. 

\bc  \label{Embedding of M(G,w) in WAP* Corollary} \rm The map $\Theta: \M(G, \om) \hookrightarrow WAP(A_\om^*)^*$, as defined in (\ref{Def of embedding map theta}), is a $so_l$-weak$^*$ and $so_r$-weak$^*$ continuous isometric homomorphic embedding that extends $\eta_{WAP}: \L^1(G, \om) \hookrightarrow WAP(A_\om^*)^*$. 
\ec 

As shown in \cite{Dal-Lau},   $WAP(A_\om^*)$ may fail to equal  $WAP(G, \om^{-1}) = \left\{f: {f \over \om} \in WAP(G)\right\}$.   Our final two results  are needed in \cite{Kro-Ste-Sto-Yee}. Corollary  \ref{Strict-Wk* continuous extension Corollary} improves \cite[Theorem 5.6]{Il-Sto} in the case of $\L^1(G, \om)$:    

\bc \label{Strict-Wk* continuous extension Corollary} \rm Let $\fB$ be a dual Banach algebra, $\vp: \L^1(G, \om) \ra \fB$ a bounded homomorphism. Then there is a unique $so_l$-weak$^*$ and $so_r$-weak$^*$ continuous homomorphic extension $\widetilde{\vp}: \M(G, \om) \ra \fB$ of $\vp$.
\ec 

\begin{proof}   Letting  $\vp_{WAP}: WAP(A_\om^*)^* \ra \fB$ be the weak$^*$-weak$^*$ continuous extension of $\vp$  from Theorem \ref{Runde Thm} and $\Theta: \M(G, \om) \hookrightarrow WAP(A_\om^*)^*$ the $so_l/so_r$-weak$^*$  continuous embedding from  Corollary \ref{Embedding of M(G,w) in WAP* Corollary}, $\widetilde{\vp} := \vp_{WAP} \circ \Theta$ is the desired extension; uniqueness follows from the $so_l$-density of $\L^1(G, \om)$ in $\M(G, \om)$.  \end{proof}

\bc \label{so-wk* conts on unit ball Cor}   \rm Let $\fB$ be a dual Banach algebra,  $\vp: \M(G, \om) \ra \fB$  a bounded homomorphism that is $so_l$-weak$^*$ continuous on the unit ball of $\M(G, \om)$. Then $\vp$ is $so_l$-weak$^*$ and $so_r$-weak$^*$ continuous on all of $\M(G, \om)$.  
 \ec 
 
 \begin{proof}   By Corollary \ref{Strict-Wk* continuous extension Corollary}, the restriction, $\vp_1$, of $\vp$ to $\L^1(G, \om)$ has  a $so_l$/$so_r$-weak$^*$  continuous  extension $\widetilde{\vp_1}: \M(G, \om) \ra \fB$.  As noted before, $\L^1(G, \om)_{\| \c \| \leq 1}$ is $so_l$-dense in $\M(G, \om)_{\| \c \| \leq 1}$, so $\vp = \widetilde{\vp_1}$ on $\M(G, \om)_{\| \c \| \leq 1}$ and therefore on $\M(G, \om)$. 
 \end{proof}  
 
 \br \rm Suppose that $(H, \om_H)$ is another weighted locally compact group and $\vp: \M(G, \om) \ra \M(H, \om_H)$ is a bounded algebra isomorphism. By \cite[Lemma 3.3]{Gha-Zad} --- which applies, as written, to $\M(G, \om)$ --- $\vp$ is $so_l$-weak$^*$ continuous on bounded subsets of $\M(G, \om)$. By Corollary \ref{so-wk* conts on unit ball Cor}, $\vp$ is $so_l/so_r$-weak$^*$ continuous on all of $\M(G, \om)$. 
  \er

\noindent {\bf Acknowledgements:} The author is grateful to Fereidoun  Ghahramani  for  helpful discussions regarding the topic of this paper.

\noindent {\sc Department of Mathematics and Statistics, University
of Winnipeg, 515 Portage Avenue, Winnipeg, MB, R3B 2E9, Canada }

\noindent email: {\tt r.stokke@uwinnipeg.ca}

\end{document}